\documentclass[11pt]{amsart}

\usepackage{amsmath,amssymb}
\usepackage{listings}

\usepackage{tikz}
\usetikzlibrary{patterns}
\usepackage[ruled,vlined,linesnumbered,norelsize]{algorithm2e}
\SetEndCharOfAlgoLine{}
\SetKwProg{Function}{Function}{}{}
\SetKwInOut{Input}{Input}
\SetKwInOut{Output}{Output}

\usepackage{hyperref} 
\usepackage{url}

\makeatletter
\newcommand{\addresseshere}{%
  \enddoc@text\let\enddoc@text\relax
}
\makeatother

\newcommand{\Z}{\mathbb{Z}}
\newcommand{\Q}{\mathbb{Q}}
\newcommand{\GFq}[1]{\mathbb{F}_{#1}}
\newcommand{\EE}{\mathcal{E}}
\newcommand{\Es}{\mathcal{E}^*}
\newcommand{\As}{A^*}
\newcommand{\Bs}{B^*}

\newcommand{\qdiff}[1]{{#1}'}

\newtheorem{lemma}{Lemma}[section]
\newtheorem{theorem}[lemma]{Theorem}
\newtheorem{proposition}[lemma]{Proposition}

\begin{document}

\title[Using CCR modular polynomials]{Using the {C}harlap-{C}oley-{R}obbins polynomials for computing isogenies}
\author{François Morain}
\address{
    LIX - Laboratoire d'informatique de l'École polytechnique 
    \emph{and}
    GRACE - Inria Saclay--Île-de-France
}
\email{morain@lix.polytechnique.fr} 
\thanks{The author is on leave from D\'{e}l\'{e}gation G\'{e}n\'{e}rale
pour l'Armement.}
\date{\today} 

\maketitle

\begin{abstract}
The SEA algorithm for computing the cardinality of elliptic curves
over finite fields in many characteristic uses modular
polynomials. These polynomials come into different flavors, and
methods to compute them flourished. Once equipped with some modular
polynomials for prime $\ell$, algebraic formulas are used to compute a
curve $\Es/\GFq{q}$ that is $\ell$-isogenous to the curve of interest
$\EE$. These formulas involve derivatives of the modular polynomial
that may sometime vanish. One way to overcome this problem is to use
alternative trivariate polynomials $U_\ell$, $V_\ell$, and $W_\ell$
introduced by Charlap, Coley and Robbins to overcome some
difficulties in the first versions of Elkies's approach. We give
properties of these polynomials, as well as formulas to compute the
isogenous curve that were sketched by Atkin. Also we investigate
another suggestion of Atkin using modular polynomials associated to a
power product of Dedekind's $\eta$ function.
\end{abstract}

\section{Introduction}

Computing isogenies is the central ingredient of the
Schoof-Elkies-Atkin (SEA) algorithm
that computes the cardinality of elliptic curves over finite fields of
large characteristic~\cite{Schoof95,Atkin92b,Elkies98} and also
\cite{BlSeSm99}. More recently, it is used in post-quantum
cryptography~\cite{ChLaGo09,JaFe11,CaLaMaPaRe18,FeKoLePeWe20} among
others, as well as the cryptosystems~\cite{Couveignes06,RoSt06,FeKiSm18}.

Given a curve $\EE/\GFq{q}$, we need to compute
$\ell$-isogenous curves for (small) prime $\ell$'s, starting from a root
of some degree $\ell+1$ modular polynomial $\Psi_\ell(X, Y)$ in $X =
j(\EE)$. The coefficients of the isogenous curve $\Es/\GFq{q}$ together
with the kernel polynomial of the isogeny are computed from
partial derivatives of $\Psi_\ell(X, Y)$. This method works except in
cases where one of the derivatives is zero. This is a theoretical as
well as practical problem. An alternative to this is to use the
original approach of Elkies \cite{Elkies98}, namely using algebraic
relations for $\As$ and $\Bs$. Another choice is to use the
CCR polynomials $(U_\ell, V_\ell, W_\ell)$ introduced in \cite{ChCoRo91},
for the price of finding the roots of three
degree $\ell+1$ polynomials instead of $1$. Several methods for
computing these polynomials, as well as rational representations of
$\As$ and $\Bs$ are given in~\cite{Morain23a} (as well as the
cited references in this article).

The aim of our work (which is closely related to ~\cite{Morain23a}) is
to give formulas to compute the necessary
parameters for Elkies's algorithms, using partial derivatives of
$U_\ell$, in the spirit of Atkin and following the suggestion in
\cite{Atkin92b}. The same work is done for the alternate polynomial
$U^a$ suggested by Atkin when $\ell \equiv 11 \bmod 12$.

The content is as follows. We review the SEA algorithm in
Section~\ref{sct:SEA}, including formulas for Eisenstein series and
introduce the CCR polynomials in
Section~\ref{sct:CCR}. Section~\ref{sct:CCRA} is the central part of
the article, working out the formulas giving the coefficients of the
isogenous curve.

\section{Schoof/Elkies/Atkin}
\label{sct:SEA}

\subsection{Prerequisites}

\subsubsection{Division polynomials}

For $\EE: y^2=x^3+A x+B$, multiplication of a point $(X, Y)$ by positive
$n$ on $\EE$ is given by
$$[n] (X, Y) = \left(\frac{\phi_n(X, Y)}{\psi_n(X, Y)^2}, \frac{\omega_n(X,
Y)}{\psi_n(X, Y)^3}\right)$$
where the polynomials satisfy
$$\phi_n = X \psi_n^2 - \psi_{n+1} \psi_{n-1},
\quad 4 Y \omega_n = \psi_{n+2} \psi_{n-1}^2 -
\psi_{n-2}\psi_{n+1}^2$$
and $\phi_n, \psi_{2n+1}, \psi_{2n}/(2Y), \omega_{2n+1}/Y,
\omega_{2n}$ belong to $\Z[A, B, X]$. It is customary to simplify this
using
$$f_n(X) = \left\{\begin{array}{ll}
\psi_n(X, Y) & \mathrm{for}\; n \;\mathrm{odd}\\
\psi_n(X, Y)/(2 Y) & \mathrm{for}\; n \;\mathrm{even}
\end{array}\right.$$
with first values
$$f_{-1} = -1, \quad f_0 = 0, \quad f_1 = 1, \quad f_2 = 1,$$
$$f_3(X, Y) = 3 X^4 + 6 A X^2 + 12 B X - A^2$$
$$f_4(X, Y) = X^6 + 5 A X^4 + 20 B X^3 - 5 A^2 X^2
- 4 A B X -8 B^2 - A^3.$$
The degree of $f_n$ is $(n^2-1)/2$ for odd $n$ and $(n^2-4)/2$ for
even $n$. If $X$ has weight $1$, $A$ weight 2 and $B$ weight 3, all
monomials in $f_n$ have the same weighted degree equal to the degree
of $f_n$.

\subsubsection{Modular functions and such}

Letting $q = \exp(2 i \pi\tau)$, one defines
$$E_2(q) = 1 -24 \sum_{n=1}^\infty \delta_1(n) q^n,$$
$$E_4(q) = 1 + 240 \sum_{n=1}^\infty \delta_3(n) q^n,$$
$$E_6(q) = 1 - 504 \sum_{n=1}^\infty \delta_5(n)q^n,$$
where $\delta_r(n)$ denotes the sum of the $r$-th powers of the
divisors of $n$. The series $E_4$ and $E_6$ are modular forms of
weight $4$ and $6$ respectively. The function $E_2$ is not a modular
form, but note that $\mathcal{F}_n(\tau)
= E_2(\tau)-n E_2(n\tau)$ is a modular form of weight 2 for
$\Gamma_0(n)$ and trivial multiplier (see 
\cite{Ramanujan14} and \cite{Berndt89}).

When $F(q) = \sum_{n\geq n_0} a_n q^n$, we introduce the operator
\begin{equation}\label{operator}
\qdiff{F}(q) = \frac{1}{2 i \pi} \, \frac{dF}{d\tau} = q \frac{d F}{d
q} = \sum_{n \geq n_0} n a_n q^n.
\end{equation}
Several identities are classical:
\begin{equation}\label{jdelta}
\Delta = \frac{E_4^3-E_6^2}{1728}, \quad
j = \frac{E_4^3}{\Delta}, \quad j - 1728 = \frac{E_6^2}{\Delta},
\end{equation}
\begin{equation}\label{formj}
\frac{\qdiff{j}}{j} = - \frac{E_6}{E_4}, \quad \frac{\qdiff{j}}{j-1728} = -
\frac{E_4^2}{E_6}, \quad \qdiff{j} = - \frac{E_4^2 E_6}{\Delta},
\quad \frac{\qdiff{\Delta}}{\Delta} = E_2,
\end{equation}
to which we add the Ramanujan differential system:
\begin{equation}\label{diff46}
{3 \qdiff{E_4}} = {E_4} E_2 -
{E_6}, \quad {2 \qdiff{E_6}} = {E_6} E_2 - {E_4^2},
\quad 12 \qdiff{E_2} = E_2^2 - E_4.
\end{equation}

\subsection{Schoof's approach}

Let $\EE/\GFq{q}$ be an elliptic curve of cardinality $m = p+1-t$ with $|t|
\leq 2\sqrt{q}$ by Hasse's theorem. Schoof gave the first
deterministic polynomial time algorithm to compute $m$. The idea is to
use the action of the Frobenius of $\EE$ on $\ell$-division points to
find $t$ via the characteristic equation $\phi^2-[t]\phi + [q] = 0$
modulo $(f_\ell, \ell)$.

\subsection{Using isogenies}

Elkies and Atkin gave subsequent improvements to make Schoof's
algorithm efficient (and probabilistic) and usable in practice.
Elkies described how to use isogenies to find factors of small degree
of $f_\ell(X)$ over a finite field, provided the Frobenius equation
splits modulo $\ell$. 
Using modular polynomials, Elkies worked out a procedure to compute all the
parameters needed to build a degree $\ell$ isogeny from $\EE: Y^2 =
X^3 + A X + B$ to some
curve $\Es: Y^2 = X^3 + \As X + \Bs$ and the kernel polynomial of the
isogeny, thereby giving
the factor we need. Atkin designed his own route towards the same
goal, putting the emphasis on the use of more modular equations for
$X_0(\ell)$ and its quotients.

One has (after renormalization):
$$A = -3 E_4(q), B = -2 E_6(q).$$
With a compatible scaling, we get
$$\As = -3 \ell^4 E_4(q^\ell), \quad \Bs = -2 \ell^6 E_6(q^\ell).$$
More importantly, writing $\sigma_r$ for the power sums of the roots of
$U_\ell$, we have
$$\sigma_1 = \frac{\ell}{2} (\ell E_2(q^\ell)-E_2(q)) =
-\frac{\ell}{2} \mathcal{F}_\ell(q).$$
Beyond this, Elkies proved that
$$A-\As = 5 (6 \sigma_2 + 2 A \sigma_0),$$
$$B-\Bs = 7 (10 \sigma_3 + 6 A \sigma_1+4 B \sigma_0),$$
together with an induction relation satisfied by other $\sigma_k$ for
$k > 3$. As a consequence $\As$ and $\Bs$ belong to $\Q[\sigma_1, A, B]$
since $\sigma_2$ and $\sigma_3$ do.

Given these quantities, there are several algorithms to get the
isogeny. We refer to \cite{BoMoSaSc08} for this.

\section{The polynomials of Charlap-Coley-Robbins}
\label{sct:CCR}

\subsection{Theory}

We start from an elliptic curve $\EE: Y^2 = X^3+A X+B$ and we fix some
odd prime $\ell$, putting $d = (\ell-1)/2$. Our aim is to find the
equation of an $\ell$-isogenous curve $\Es: Y^2 = X^3 + \As X + \Bs$. 

\begin{theorem}
There exist three polynomials $U_\ell$, $V_\ell$,
$W_\ell$ in $\Z[X, Y, Z, 1/\ell]$ of degree $\ell+1$ in $X$
such that $U_\ell(\sigma_1, A, B)=0$, respectively $V_\ell(\As, A, B) = 0$,
$W_\ell(\Bs, A, B) = 0$.
\end{theorem}

Let us turn our attention to the properties of these polynomials.
\begin{theorem}
When $\ell > 3$, the polynomials $U_\ell$, $V_\ell$, $W_\ell$ live in
$\Z[X, Y, Z]$.
\end{theorem}

\begin{proposition}
Assigning respective weights 1, 2, 3 to $X$, $Y$, $Z$,
the monomials in $U_\ell$, $V_\ell$ and $W_\ell$ have generalized
degree $\ell+1$.
\end{proposition}

\subsection{Computing isogenous curves over finite fields}
\label{sct:summary}

When using $U_\ell$, $V_\ell$, $W_\ell$, we need to find the roots of
three polynomials of degree $\ell+1$ instead of $1$. In general, if
$U_\ell$ has rational roots (it should be 1, $2$ or $\ell+1$), then
this is the case for each of $V_\ell$, $W_\ell$. For each triplet of
solutions $(\sigma_1, z_1, z_2)$ we need to test whether this leads to
an isogeny or not. To speed up things, we may compute rational
fractions for $\As$ and $\Bs$ as explained in \cite{NoYaYo20} (see
also \cite[\S 7]{PoSc13}). Another path was sketched by Atkin in
\cite{Atkin92b}, and this is what we describe next.

\section{Revisiting CCR {\em \`a la} Atkin}
\label{sct:CCRA}

The idea is to generalize the approach in
\cite{Atkin88b,Atkin92b,Morain95a}, that is exploit $q$-series
identities to get the parameters $(\sigma, \As, \Bs)$, where we write
$\sigma$ for $\sigma_1$ from now on.

\subsection{Properties of $U_\ell$}

We write for readability $U = U_\ell$ and
$$\partial_\sigma = \frac{\partial U}{\partial \sigma},
\partial_4 = \frac{\partial U}{\partial E_4},
\partial_6 =  \frac{\partial U}{\partial E_6}.$$
and propagate the notation to double derivatives.

The polynomial $U$ is homogeneous with weights, so that
\begin{equation}\label{homU}
(\ell+1) U = \sigma \partial_\sigma + 2 E_4 \partial_4 + 3 E_6
\partial_6.
\end{equation}
Note that partial derivatives of $U$ are also homogeneous polynomials
and we find
\begin{eqnarray}\label{dxx}
\ell \partial_\sigma &=& \sigma \partial_{\sigma\sigma} + 2 E_4
\partial_{\sigma 4} + 3 E_6 \partial_{\sigma 6}, \\
(\ell-1) \partial_4 &=& \sigma \partial_{\sigma 4} + 2 E_4
\partial_{44} + 3 E_6 \partial_{46}, \\
(\ell-2) \partial_6 &=& \sigma \partial_{\sigma 6} + 2 E_4
\partial_{46} + 3 E_6 \partial_{66}.
\end{eqnarray}

\subsection{Getting the isogenous curve from CCR polynomials}

\subsubsection{Finding $\tilde{E}_4$}

\begin{proposition}
The value of $\tilde{E}_4$ is given by
$$- \frac{4 \ell (3 E_4^2 \partial_6+2E_6 \partial_4)
- \partial_\sigma (\ell^2 E_4+4  \sigma^2)}
{\ell^4 \partial_\sigma}.$$
\end{proposition}

\medskip
\noindent
{\em Proof:}
We differentiate (using (\ref{operator})) $U(\sigma,
E_4, E_6) = 0$ to get
\begin{equation}\label{eq1}
\sigma' \partial_\sigma + E_4' \partial_4 + E_6' \partial_6 = 0.
\end{equation}
We differentiate $\sigma = \ell/2 (\ell \tilde{E}_2 - E_2)$
leading to
$$\sigma' = \frac{\ell}{2} \; (\ell^2 \tilde{E}_2'- E_2') =
\frac{\ell}{24}\; (\ell^2 (\tilde{E}_2^2-\tilde{E}_4) -
(E_2^2-E_4)).$$
Replace $\ell \tilde{E}_2$ by $2\sigma/\ell+ E_2$ to get
$$\sigma' = \frac{\ell}{24} \; \left(\frac{4 \sigma^2}{\ell^2} +
\frac{4 \sigma}{\ell} E_2 - (\ell^2 \tilde{E}_4- E_4)\right),$$
that we plug in (\ref{eq1}) together with the expressions for
$\qdiff{E_4}$ and $\qdiff{E_6}$ from equation (\ref{diff46}) to get a
polynomial of degree 1 in $E_2$ whose coefficient of $E_2$ is
$$\sigma \partial_\sigma + 2 E_4 \partial_4 + 3 E_6 \partial_6,$$
which we recognize in (\ref{homU}). Therefore, we get
\begin{equation}\label{eq2}
(\ell+1) U E_2 + \frac{\ell\,\partial_\sigma}{4} \left(\frac{4
\sigma^2}{\ell^2} - (\ell^2 \tilde{E}_4- E_4)\right)  - 2
E_6 \partial_4 - 3 E_4^2 \partial_6 = 0
\end{equation}
from which we deduce $\tilde{E}_4$ since $U(\sigma, E_4, E_6) = 0$. $\Box$

\subsubsection{Finding $\tilde{E}_6$}

\begin{proposition}
The value of $\tilde{E}_6$ may be written
$$\tilde{E}_6 = -\, \frac{N}{\ell^6 \, \partial_\sigma^3}$$
where $N$ is some polynomial of degree 3 in $\ell$ and given at the
end of the proof.
\end{proposition}

\medskip
\noindent
{\em Proof:}
We differentiate (\ref{eq1}).
\begin{eqnarray}
\label{ea1}
& & \sigma'' \partial_\sigma + {\sigma'}(\sigma' \partial_{\sigma\sigma}
+ E_4' \partial_{\sigma 4} + E_6' \partial_{\sigma 6}) \\ 
\label{ea2} &+&E_4'' \partial_4 + E_4' (\sigma' \partial_{4\sigma} + E_4'
\partial_{44} + E_6' \partial_{46}) \\
 \label{ea3}&+& E_6'' \partial_6 + E_6' (\sigma' \partial_{6\sigma} + E_4'
\partial_{64} + E_6' \partial_{66}) = 0
\end{eqnarray}
We compute in sequence
$$12 E_2'' = 2 E_2 E_2' - E_4' = E_2 (E_2^2-E_4)/6 - (E_2 E_4 - E_6)/3,$$
$$12 {\tilde{E}_2}'' = 2 \tilde{E}_2 \tilde{E}_2' - \tilde{E}_4'
= \tilde{E}_2 (\tilde{E}_2^2-\tilde{E}_4)/6 - (\tilde{E}_2
\tilde{E}_4 - \tilde{E}_6)/3,$$
which give us the value
$$\sigma'' = \frac{\ell}{2} \; (\ell^3 {\tilde{E}_2}''- E_2'')$$
to be used in (\ref{ea1}).
Differentiating relations of (\ref{diff46}), we get
$$E_4'' = \frac{1}{3} \, (E_2' E_4 + E_2 E_4' - E_6'),
\quad E_6'' = \frac{1}{2} \, (E_2' E_6 + E_2 E_6' - 2 E_4 E_4'),$$
to be used in lines (\ref{ea2}) and (\ref{ea3}) respectively. We
replace $\tilde{E}_4$ by its value from (\ref{eq2}), and $\tilde{E}_2$ using
$\sigma = (\ell/2) ( \ell \tilde{E}_2 - E_2)$.
This finally yields an expression
as polynomial in $E_2$: $$C_2 E_2^2 + C_1 E_2 + C_0 = 0.$$
The unknown $\tilde{E}_6$ is to be found in $C_0$ only.

By luck(?)
\begin{proposition}
The coefficients $C_1$ and $C_2$ vanish for a triplet such that
$U_\ell(\sigma, E_4, E_6) = 0$.
\end{proposition}

\medskip
\noindent
{\em Sketch of the proof:} The strategy to
prove this is the same in both cases. Replace
$\partial_{\sigma\sigma}$, $\partial_{44}$ and $\partial_{66}$ by
their values from (\ref{dxx}). Factoring the resulting expressions
yields the same factor $\sigma \partial_\sigma + 2 E_4 \partial_4 + 3
E_6 \partial_6$, which cancels $C_1$ and $C_2$. We add a \verb+SageMath+
script for the convenience of the reader as an appendix to this
work. $\Box$

We are left with 
$$\tilde{E}_6 = -\, \frac{N}{\ell^6 \, \partial_\sigma^3}$$
where $N$ is a polynomial in degree 3 in $\ell$
$$N = -E_6 \partial_\sigma^3 \ell^3 + c_2 \ell^2 + 12
\partial_\sigma^2 \sigma (3 E_4^2 \partial_6+2 E_6 \partial_4)
\ell - \partial_\sigma^3 \sigma^3.$$
The coefficient $c_2$ is heavy looking and we give slightly factored
as a polynomial in $E_4$:
\begin{eqnarray*}
c_2 &=&
18 (\partial_6^2 \partial_{\sigma\sigma}-2 \partial_6 \partial_\sigma \partial_{\sigma 6}+ \partial_{66} \partial_\sigma^2) E_4^4\\
&&+(24 E_6 \partial_4 (\partial_6 \partial_{\sigma\sigma}-
\partial_\sigma \partial_{\sigma 6})
+24 E_6 \partial_\sigma (\partial_{46} \partial_\sigma- \partial_6 \partial_{\sigma 4})+10 \partial_4 \partial_\sigma^2) E_4^2\\
&&+3 \partial_\sigma^2 (7 E_6 \partial_6 - \sigma \partial_\sigma) E_4
+8 E_6^2 (\partial_4^2 \partial_{\sigma\sigma}-2 \partial_4
\partial_\sigma \partial_{\sigma 4}+\partial_{44} \partial_\sigma^2). \Box
\end{eqnarray*}

\subsubsection{Numerical example}

Consider $E: Y^2 = X^3+X+3$ over $\GFq{1009}$ and $\ell = 5$. Using
$$U_5(X) = X^6 + 20 X^4 A + 160 X^3 B - 80 X^2 A^2 - 128 X A B - 80
B^2,$$
we select $\sigma = 584$ and compute
$$\partial_\sigma=905, \partial_4=779, \partial_6=140$$
from which $\tilde{E}_4=497$, $\As = 441$. After tedious computations,
we find $\Bs = 997$.

\subsection{The case $\ell \equiv 11 \bmod 12$}

In this case, Atkin suggests to replace $\sigma$ with $f(q) = 
(\eta(q)\eta(q^\ell))^2$, where $\eta$ is Dedekind's function. The
corresponding modular polynomial
$U_\ell^a$ can be computed using the techniques described in \cite{Morain23a}.
For instance (using the basis with $E_4$, $E_6$ and $\Delta$):
$$U_{11}^a(X)=
{X}^{12}-990\,\Delta\,{X}^{6}+440\,E_4\,\Delta\,{X}^{4}-165\,E_6\,\Delta\,{X}^{3}+22\,E_4^{2}\,\Delta{X}^{2}$$
$$-E_6\,E_4\,\Delta\,X-11\,{\Delta}^{2},$$
which is sparser than $U_{11}(X)$.

\subsubsection{Some properties of $U_\ell^a$}

Let $v_2$ (resp. $v_3$) be the maximal power of $1/2$ (resp. $1/3$) of
the coefficients of $U_\ell^a(X)$ that we found experimentally
$$\begin{array}{|c|c|c|} \hline
\ell & v_2 & v_3 \\ \hline
11 & 16 & 12 \\
23 & 32 & 24 \\
47 & 64 & 48 \\
59 & 80 & 60 \\
71 & 96 & 72 \\
\hline
\end{array}$$
It seems that $f/12$ should be a more sensible choice, leading to
$12^{\ell+1} U_\ell^a(X/12)$ having integer coefficients.

We have formulas analogous to (\ref{dxx}), due the corresponding
homogeneous property
\begin{eqnarray}\label{dff}
\ell \partial_f &=& f \partial_{ff} + 2 E_4
\partial_{f 4} + 3 E_6 \partial_{f 6}, \\
(\ell-1) \partial_4 &=& f \partial_{f 4} + 2 E_4
\partial_{44} + 3 E_6 \partial_{46}, \\
(\ell-2) \partial_6 &=& f \partial_{f 6} + 2 E_4
\partial_{46} + 3 E_6 \partial_{66}.
\end{eqnarray}

\subsubsection{Computing $\sigma$, $\tilde{E}_4$ and $\tilde{E}_6$}

\begin{proposition}
The value of $\sigma$ is
$$\sigma = {\frac {{\ell}\, \left( 3\,{ \partial_6}\,E_4^{2}+2\,{
\partial_4}\, E_6 \right) }{f{ \partial_f}}}.$$
\end{proposition}

\medskip
\noindent
{\em Proof:}
Remark that $f^{12} = \Delta(z) \Delta(\ell z)$ and therefore
we deduce the discriminant $\tilde{\Delta} = f^{12}/\Delta$ of the
isogenous curve. We have also (using (\ref{formj})):
$$12 \, \frac{f'}{f} = \frac{\Delta'}{\Delta} + \ell
\frac{\tilde{\Delta}'}{\tilde{\Delta}} = E_2 + \ell \tilde{E}_2,$$
from which we deduce $f'$. Again, $U_\ell^a$ is homogeneous with
weight $\ell+1$, so that we have identities similar to those for
$\sigma$. In particular
\begin{equation}\label{homUA}
(\ell+1) U_\ell^a = f \partial_f + 2 E_4 \partial_4 + 3 E_6
\partial_6.
\end{equation}
Starting from $f' \partial_f + E_4' \partial_4 + E_6' \partial_6 = 0$,
and replacing by the known values, we find
$$(f \partial_f + 4 E_4 \partial_4 + 6 E_6 \partial_6) \, E_2
+ f \ell \tilde{E}_2 \partial_f - 6 E_4^2 \partial_6 -4 E_6 \partial_4
= 0,$$
which is
$$ f \partial_f (\ell \tilde{E}_2 - E_2) - 6 E_4^2
\partial_6 -4 E_6 \partial_4 = 0,$$
and this gives us the result. $\Box$

\begin{proposition}
The value of $\tilde{E}_4$ is given by
$$\tilde{E}_4 = - \frac{M}{\ell^2 f^2 E_4 E_6 \partial_f^3 }$$
where $M$ is a polynomial given at the end of the proof.
\end{proposition}

\medskip
\noindent
{\em Proof:}
We differentiate $f'$ to obtain:
\begin{eqnarray*}
f'' &=& \frac{1}{12} \left(f' (\ell \tilde{E}_2 + E_2)
+f (\ell^2 \tilde{E}_2' +
E_2')\right)\\
&=& \frac{f}{12^2} \left((\ell \tilde{E}_2 + E_2)^2+
\ell^2 (\tilde{E}_2^2-\tilde{E}_4) + (E_2^2-E_4))\right).
\end{eqnarray*}
We inject this together and the diagonal derivatives of (\ref{dff})
and $\tilde{E}_2 = (E_2 + 2 \sigma/\ell)/\ell$ into
\begin{eqnarray}
\label{eaa1}
& & f'' \partial_f + {f'}(f' \partial_{ff}
+ E_4' \partial_{f 4} + E_6' \partial_{f 6}) \\ 
\label{eaa2} &+&E_4'' \partial_4 + E_4' (f' \partial_{4f} + E_4'
\partial_{44} + E_6' \partial_{46}) \\
 \label{eaa3}&+& E_6'' \partial_6 + E_6' (f' \partial_{6f} + E_4'
\partial_{64} + E_6' \partial_{66}) = 0
\end{eqnarray}
to get a polynomial of degree 2 in $E_2$ whose coefficients of
degree 2 and 1 turn out to vanish. We are left with
$$\tilde{E}_4 = - \frac{M}{\ell^2 f^2 E_4 E_6 \partial_f^3 }$$
where
\begin{eqnarray*}
M &=&
24 (3 E_6 \partial_6^2 \partial_{f4}+\partial_{46} \partial_f^2 f) E_4^6\\
&&+12 (9 E_6^2 \partial_6^2 \partial_{f6}-3 E_6 \partial_6^2
\partial_f \ell+6 E_6 \partial_6 \partial_f \partial_{f6} f-\partial_6
\partial_f^2 \ell f+\partial_f^2 \partial_{f6} f^2 \\
&& \qquad -6 E_6 \partial_6^2 \partial_f+2 \partial_6 \partial_f^2 f) E_4^5
+96 E_4^4 
E_6^2 \partial_4 \partial_6 \partial_{f4}\\
&&+4 E_6 (36 E_6^2 \partial_4 \partial_6 \partial_{f6}-12 E_6
\partial_4 \partial_6 \partial_f \ell+12 E_6 \partial_4 \partial_f
\partial_{f6} f -12 E_6 \partial_{46} \partial_f^2 f\\
 &&\qquad +12 E_6 \partial_6 \partial_f \partial_{f4} f-24 E_6
 \partial_4 \partial_6 \partial_f-5 \partial_4 \partial_f^2 f) E_4^3\\
&& +E_6 (32 E_6^2 
\partial_4^2 \partial_{f4}-42 E_6 \partial_6 \partial_f^2
 f+\partial_f^3 f^2) E_4^2\\
&&+16 E_6^3 \partial_4 (3 E_6 \partial_4 \partial_{f6}-\partial_4 \partial_f \ell+2 \partial_f
 \partial_{f4} f-2 \partial_4 \partial_f) E_4\\
&&+24 E_6^4 \partial_{46} f \partial_f^2-8 E_6^3 \partial_4 \ell f \partial_f^2+8 E_6^3 \partial_{f4} f^2 \partial_f^2
+8 E_6^3 \partial_4 f \partial_f^2. \quad \Box
\end{eqnarray*}

Finally, we remark that $\Bs$ satisfies
\begin{equation}\label{twopol}
{\Bs}^2 + 6912 \tilde{\Delta}-4 {\As}^3/27 = 0,
\quad U_\ell^a(-\ell f, \As, \Bs) = 0,
\end{equation}
the latter relation coming from applying the Atkin-Lehner involution
to the modular form for $\sigma_1$.
The gcd of these two polynomials should reveal
$\Bs$. In the rare case where this gcd has degree 2 (which would imply
two elliptic curves being isogenous to $E$), we would
be forced to use higher differentials, which would look like a
formidable task.

\subsubsection{Numerical example}

Consider $E: Y^2 = X^3+X+3$ over $\GFq{1009}$. The polynomial
$U_{11}^a$ has two roots: $65$ and $333$. We take $f = 65$. We first
compute $\sigma=75$. Then $\tilde{E}_4=532$. The gcd of the two
polynomials in (\ref{twopol}) has degree 1 and root $\Bs=460$.

\section{Conclusions}

We have completed the task suggested by Atkin for using the CCR
polynomials in building isogenies. All these formulas require
$O(\ell^2)$ multiplications in the base field, due to the computation
of partial derivatives of polynomials of degree $O(\ell)$. Note that
this is the same cost as using the rational fractions giving $\As$ and
$\Bs$, but less storage is needed.

As a consequence, we have several algorithms and
formulas to be used, depending on the practical problem to be solved.

\def\noopsort#1{}\ifx\bibfrench\undefined\def\biling#1#2{#1}\else\def\biling#1#2{#2}\fi\def\Inpreparation{\biling{In
  preparation}{en
  pr{\'e}paration}}\def\Preprint{\biling{Preprint}{pr{\'e}version}}\def\Draft{\biling{Draft}{Manuscrit}}\def\Toappear{\biling{To
  appear}{\`A para\^\i tre}}\def\Inpress{\biling{In press}{Sous
  presse}}\def\Seealso{\biling{See also}{Voir
  {\'e}galement}}\def\Editor{\biling{Ed.}{R{\'e}d.}}

\appendix

\section{A script to check the computations}

This \verb+SageMath+~\cite{sagemath} script can also be downloaded
from the author's web page.

{\def\baselinestretch{0.9}
\begin{lstlisting}[language=python]
# This script is devoted to the computation and verification of several
# identities related to CCR polynomials using the notations of the preprint.

# one ring to rule them all
R.<ell,E2,E4,E6,sigma,E4t,E6t,d4,d6,s,ds,ds4,ds6,d46,f,df,df4,df6>
    =PolynomialRing(Rationals(),18)

########## The CCR case

# returns ell^-4 * ds^-1 * (-12*ell*E4^2*d6 + ...)
def check_E4t():
    E4p=(E2*E4 - E6)/3
    E6p=(E2*E6-E4^2)/2
    E2p=(E2^2-E4)/12
    E2t=(E2+2*sigma/ell)/ell
    sigp=ell/24*(4*sigma^2/ell^2+4*sigma/ell*E2-(ell^2*E4t-E4))
    tmp=sigp*ds+E4p*d4+E6p*d6
    tmp=tmp.numerator()
    print("degree(tmp, E2)=", tmp.degree(E2))
    # check that coeff of E2 is zero
    c1=tmp.coefficient({E2:1})
    # is a multiple of (2*E4*d4 + 3*E6*d6 + f*df), hence 0
    print("c1=", c1.factor())
    # find sigma as a root of constant coefficient
    e4t=tmp.coefficient({E2:0})
    e4t=-e4t.coefficient({E4t:0})/e4t.coefficient({E4t:1})
    # sig contains the value of sigma
    return e4t.factor()

# returns
# ell^-6 * ds^-3 * sigma^-1 * E6^-1 * E4^-1 * (-18*ell^3*E4^5*E6*d6^2*ds+...)
def check_E6t():
    e4t=check_E4t()
    E4p=(E2*E4 - E6)/3
    E6p=(E2*E6-E4^2)/2
    E2p=(E2^2-E4)/12
    E2t=(E2+2*sigma/ell)/ell
    sigp=ell*(4*sigma^2/ell^2+4*sigma/ell*E2-(ell^2*e4t-E4))/24
    # more derivatives
    E4pp=1/3*(E2p*E4+E2*E4p-E6p)
    E6pp=1/2*(E2p*E6+E2*E6p-2*E4*E4p)
    # crucial values
    E4tp=1/3*(E2t*e4t-E6t)
    E2tp=(E2t^2-e4t)/12
    E2pp=1/12*(2*E2*E2p-E4p)
    E2tpp=1/12*(2*E2t*E2tp-E4tp)
    sigpp=ell*(ell^3*E2tpp-E2pp)/2
    # inject diagonal derivatives
    dss = (ell*ds    -2*E4*ds4 -3*E6*ds6)/sigma
    d44 = ((ell-1)*d4-sigma*ds4-3*E6*d46)/(2*E4)
    d66 = ((ell-2)*d6-sigma*ds6-2*E4*d46)/(3*E6)
    # starting point
    tmp=      sigpp*ds+sigp*(sigp*dss+E4p*ds4+E6p*ds6)
    tmp=tmp + E4pp*d4+E4p*(sigp*ds4+E4p*d44+E6p*d46)
    tmp=tmp + E6pp*d6+E6p*(sigp*ds6+E4p*d46+E6p*d66)
    tmp=tmp.numerator()
    c2=tmp.coefficient({E2:2})
    print("E6t.c2=", c2.factor())
    c1=tmp.coefficient({E2:1})
    print("E6t.c1=", c1)
    c0=tmp.coefficient({E2:0})
    e6t=-c0.coefficient({E6t:0})/c0.coefficient({E6t:1})
    return e6t.factor()

########## The case of ell = 11 mod 12, Atkin's variant
def check11_sigma():
#    R.<ell,E2,E4,E6,sigma,d4,d6,f,df>=PolynomialRing(Rationals(),9)
    tmp=2*E4*d4+3*E6*d6+f*df
    E4p=(E2*E4 - E6)/3
    E6p=(E2*E6-E4^2)/2
    E2p=(E2^2-E4)/12
    E2t=(E2+2*sigma/ell)/ell
    fp=f/12*(ell*E2t+E2)
    tmp=fp*df+E4p*d4+E6p*d6
    tmp=tmp.numerator()
    # check that coeff of E2 is zero
    c1=tmp.coefficient({E2:1})
    # is a multiple of (2*E4*d4 + 3*E6*d6 + f*df), hence 0
    print("c1=", c1.factor())
    # find sigma as a root of constant coefficient
    sig=tmp.coefficient({E2:0})
    sig=-sig.coefficient({sigma:0})/sig.coefficient({sigma:1})
    # sig contains the value of sigma
    return sig.factor()

# returns
# (-1) * df^-3 * f^-2 * ell^-2 * E6^-1 * E4^-1 * (-36*ell*E4^5*E6*d6^2*df+...)
def check11_E4t():
    sig=check11_sigma()
    E4p=(E2*E4 - E6)/3
    E6p=(E2*E6-E4^2)/2
    E2p=(E2^2-E4)/12
    E2t=(E2+2*sig/ell)/ell
    fp=f/12*(ell*E2t+E2)
    fpp=f/12^2*((ell*E2t+E2)^2+ell^2*(E2t^2-E4t)+(E2^2-E4))
    E4pp=1/3*(E2p*E4+E2*E4p-E6p)
    E6pp=1/2*(E2p*E6+E2*E6p-2*E4*E4p)
    # inject diagonal derivatives
    dff = (    ell*df-2*E4*df4 -3*E6*df6)/f
    d44 = ((ell-1)*d4-f*df4-3*E6*d46)/(2*E4)
    d66 = ((ell-2)*d6-f*df6-2*E4*d46)/(3*E6)
    tmp=       fpp*df+ fp*(fp*dff+E4p*df4+E6p*df6)
    tmp=tmp + E4pp*d4+E4p*(fp*df4+E4p*d44+E6p*d46)
    tmp=tmp + E6pp*d6+E6p*(fp*df6+E4p*d46+E6p*d66)
    tmp=tmp.numerator()
    print("degree(tmp, E2)=", tmp.degree(E2))
    c2=tmp.coefficient({E2:2})
    print("E4t.c2=", c2.factor())
    c1=tmp.coefficient({E2:1})
    print("E4t.c1=", c1)
    c0=tmp.coefficient({E2:0})
    e4t=-c0.coefficient({E4t:0})/c0.coefficient({E4t:1})
    return e4t.factor()
\end{lstlisting}}

\end{document}